\documentstyle[remark,amssymb,12pt]{article}
\newbox\smilebox
\newbox\anchorbox
\newbox\noanchorbox
\newbox\tempbox

\setbox\smilebox=\hbox{$\smile$}

\def\anchor{\hbox{\vtop{
           \hbox to \wd\smilebox{\hfil\vrule width.4pt height7pt depth1pt\hfil}
           \vskip  -11.5truept
           \hbox to \wd\smilebox{\hfil$\smile$\hfil}}}}
\setbox\anchorbox=\anchor
\def\noanchor{\hbox{\vtop{
           \hbox to \wd\anchorbox{\hfil\anchor\hfil}
           \vskip -14truept
           \hbox to \wd\anchorbox{\hfil/\hfil}}}}
\setbox\noanchorbox=\noanchor

\def\fg#1#2#3{\setbox\tempbox=\hbox{$\scriptstyle{#2}$}
\ifnum\wd\anchorbox>\wd\tempbox\dimen255=\wd\anchorbox
\else\dimen255=\wd\tempbox\fi
{#1\,\vtop{\hbox to \dimen255{\hfil\anchor\hfil}
           \vskip -6truept
           \hbox to \dimen255{\hfil$\scriptstyle{#2}$\hfil}}
           \,#3}}

\def\nfg#1#2#3{\setbox\tempbox=\hbox{$\scriptstyle{#2}$}
\ifnum\wd\noanchorbox>\wd\tempbox\dimen255=\wd\noanchorbox
\else\dimen255=\wd\tempbox\fi
{#1\,\vtop{\hbox to \dimen255{\hfil\noanchor\hfil}
           \vskip -6truept
           \hbox to \dimen255{\hfil$\scriptstyle{#2}$\hfil}}
           \,#3}}

\setbox1=\hbox{$\bot$}

\def\north#1#2{#1\,
\hbox{$\bot$\llap {\hbox to\wd1 {\hfil $/$\hfil}}}
\,#2}

\def\nao#1#2#3{#1\  \hbox{\vtop{ 
\baselineskip=4pt
\hbox{$\bot$\llap {\hbox to\wd1 {\hfil $/$\hfil}}
\hskip .05em \llap{\hbox{$^{\scriptscriptstyle{a}}$}}}\hbox{$\scriptstyle
{#2}$}}}\, #3}

\def\bp{\par{\bf Proof.}$\ \ $}
\def\Claim#1{\par\bigskip{\bf #1.}$\ \ $}

\def\includeE#1{{\lhook\kern-3.5pt\joinrel\smash{
    \mathop{\longrightarrow}\limits^{#1}}}}

\def\efor/{Example~\ref{E4}}

\def\ep{\par\bigskip}

\def\BL/{Baldwin--Lachlan}
\def\Bu/{Buechler}
\def\Hr/{Hrushovski}
\def\lm/{locally modular}
\def\wm/{weakly minimal}
\def\nm/{non--modular}
\def\tt/{totally transcendental}
\def\ss/{superstable}
\def\ud/{unidimensional}
\def\sm/{strongly minimal}

\def\abar{\overline{a}}

\def\bbar{\overline{b}}
\def\cbar{\overline{c}}
\def\dbar{\overline{d}}

\def\hbar{\overline{h}}

\def\vbar{\overline{v}}

\def\xbar{\overline{x}}

\def\zbar{\overline{z}}

\def\acl{{\rm acl}}

\def\dom{{\rm dom}}

\def\tp{{\rm tp}}
\def\stp{{\rm stp}}

\def\tr/{trivial}
\def\nt/{non--trivial}
\def\st/{strong type}

\def\TV/{Tarski--Vaught}

\def\sc/{sound construction}
\def\ac/{atomic construction}
\def\Ceq{\bigc^{\rm eq}}

\def\fal/{functional}
\def\upl/{unique parallel lines}
\def\chp/{categorical in a higher power}

\def\text#1{\ \hbox{#1}\ }

\def\sbq{\subset}
\def\contains{\supseteq}
\def\<{\langle}
\def\>{\rangle}
\def\conc{\widehat{~~}}
\def\forces{\mathrel{\raise.4ex\hbox{$\scriptstyle \vert$\hskip-.5ex}\vdash}}

\def\K/{${\cal K}$}
\def\KM{{\cal K}}

\def\r{\restriction}
\def\rl{\r L}
\def\KMOD/{$\KM={\rm Mod}(T_1)\rl$}
\def\LK/{$L(\KM)$}
\def\S{{\rm Spec}(\KM)}
\def\al{\alpha}
\def\b{\beta}
\def\d{\delta}

\def\e{\emptyset}
\def\g{\gamma}

\def\l{\lambda}
\def\o{\omega}

\def\G{\Gamma}
\def\SN/{\hbox{$\S\ne\e$}}
\def\ao{\aleph_0}
\def\a1{\aleph_1}

\def\PCD/{${\rm PC}_{\a0}$}
\def\PC/{${\rm PC}_{\Delta}$}
\def\ra{\rightarrow}

\def\bg{\beth_\g}

\def\ss{\smallsetminus}

\def\SG/{\hbox{$\S\cap\bg$}}
\def\SNG/{\hbox{$\S\cap\bg\not =\e$}}

\def\A{{\cal A}}
\def\B{{\cal B}}
\def\C{{\cal C}}
\def\D{{\cal D}}

\def\P{{\cal P}}

\def\endproof{\enspace\vrule height6pt width4pt depth0pt\ep} 
\def\sbq{\subseteq}

\def\Khom{{\bf K}_{{\rm hom}}}
\def\P{{\cal P}}
\def\Q{{\cal Q}}
\def\R{{\cal R}}
\def\V{{\cal V}}
\def\QQ{{\Bbb Q}}
\def\RR{{\Bbb R}}
\def\ZZ{{\Bbb Z}}

\def\S{{\cal S}}
\def\A{{\frak A}}
\def\B{{\frak B}}
\def\C{{\frak C}}
\def\Ceq{\C^{\rm eq}}

\def\f{\tilde}
\def\ff{\f{f}}

\newtheorem{theorem}{Theorem}[section]

\newtheorem{lemma}[theorem]{Lemma}

\newremark{defn}[theorem]{Definition}
\newremark{example}[theorem]{Example}
\newremark{remark}[theorem]{Remark}
\author{
M. C. Laskowski\thanks{\hbox{Partially supported by an NSF 
Postdoctoral Fellowship
and Research~Grant~DMS~9403701.}}
\\Department of Mathematics\\
University of Maryland
\and
S. Shelah\\
Department of Mathematics\\
Hebrew University of Jerusalem
\and Department of Mathematics\\
Rutgers University
\thanks{The authors thank the U.S.-Israel
Binational Science Foundation for its support of this project.
This is item 518 in Shelah's bibliography.}}
\date{\today}
\title{Forcing Isomorphism II}

\begin{document}

\maketitle 

\begin{abstract}
If $T$ has only countably many complete types, 
yet has a type of infinite multiplicity
then there is a c.c.c.\ forcing notion $\Q$ such that, 
in any $\Q$-generic extension
of the universe, there are non-isomorphic models $M_1$ and $M_2$ of $T$
that can be forced isomorphic by a c.c.c.\ forcing.
We give examples showing that the hypothesis on the number of complete types is
necessary and what happens if `c.c.c.' is replaced other cardinal-preserving
adjectives.
We also give an example showing that membership in a pseudo-elementary
class can be altered by very simple cardinal-preserving forcings.

\end{abstract}

\section{Introduction}
The fact that the isomorphism type of models of a theory
can be altered by forcing was first noted by Barwise in~\cite{Barwise}.
He observed that the natural back-and-forth system obtained from
a pair of $L_{\infty,\o}$-equivalent structures gives rise to a
partial order that makes the structures isomorphic in 
any generic extension of the universe.  
Restricting attention to partial orders with desirable
combinatorial properties, e.g., the countable chain condition (c.c.c.)
and asking which theories have a pair of non-isomorphic models
that can be forced isomorphic by such a forcing
provides us with an 
alternate approach
to a fundamental question of model theory.  
The question, roughly stated, asks
which (countable) theories admit a `structure theorem' for
the class of models of the theory?
Part of the research on this question has been to discover
a definition of the phrase `structure theorem' that leads to
a natural dichotomy between theories.
In 
\cite{Shelahbook2nd} 
the second author succeeds in
characterizing the theories that have the maximal number of 
non-isomorphic models 
in every uncountable cardinality and is near a characterization
of the theories which have families of $2^\kappa$ 
pairwise non-embeddable models of size $\kappa$.
These abstract results imply the impossibility of structure theorems
(for virtually every
definition of `structure theorem') by the sheer size and complexity of
the class of models of such a theory.
On the positive side, he defines a {\it classifiable\/} theory  
(i.e., superstable, without the Dimensional Order Property (DOP)
and without the Omitting Types Order Property (OTOP))
and shows that any model of a classifiable theory can be
described in terms of an independent  tree of countable 
elementary substructures.  That is, the class of models of 
such a theory has a structure theorem in a certain sense.
In \cite{Shelah220} he analyzes which structures (of a 
fixed cardinality) can be
determined up to isomorphism by their 
Scott sentences in various infinitary languages 
(e.g., $L_{\infty,\kappa}$).

In both this paper and in
\cite{BaldLasSheforceiso}, 
we concentrate on  systems of 
invariants that are preserved under c.c.c.\ forcings
and ask which theories have their models described up to isomorphism
by invariants of this sort.
It is well-known that c.c.c.\ forcings preserve 
cardinality and cofinality, yet such forcings typically
add new subsets of
$\o$ (reals) to the universe.
We call two structures 
{\it potentially isomorphic\/}
if they can be forced isomorphic by a forcing with the countable chain
condition (c.c.c.).
The relevance of this notion is that
the existence of a pair of non-isomorphic, potentially isomorphic
structures within a class {\bf K} 
(either in the ground universe or
in a c.c.c.\ forcing extension)
implies that the isomorphism type of elements of {\bf K}
cannot be described by a c.c.c.-invariant
system of invariants.

In \cite{BaldLasSheforceiso}, 
it was shown that for countable theories $T$, if $T$ is not classifiable
then there are
non-isomorphic, potentially isomorphic models of $T$ of size $2^\o$.
In addition, certain classifiable theories were shown to have such a pair of
models.
The main theorem of this paper, Theorem \ref{big}, states that if 
$T$ is superstable,
$D(T)$ is countable (i.e., $T$ has at most $\aleph_0$ $n$-types for each $n$),
but has a type of infinite multiplicity (equivalently, $T$ is not 
$\aleph_0$-stable) then there is a c.c.c.\ forcing $\Q$
such that $\forces_\Q$ ``There are two non-isomorphic, potentially isomorphic
models of $T$.''
Combining this with the results from 
\cite{BaldLasSheforceiso}  
yields the theorem mentioned in the abstract.
We remark that the more natural question of whether any such
theory has a pair of
non-isomorphic, potentially isomorphic models
in the ground universe (as opposed to in a forcing extension)
remains open.

A consequence of these results is that the system of invariants
for the isomorphism type of 
a model of a classifiable theory mentioned above cannot be 
simplified significantly.  In particular, we conclude that
if $T$ is classifiable but not $\ao$-stable and
if  $D(T)$ is countable then the models of $T$ cannot be 
described by independent trees of {\it finite\/} subsets, for any
such tree would be preserved by a c.c.c.\ forcing.

The idea of the proof of Theorem~\ref{big}
is to build two models of $T$, each realizing
a suitably generic subset of the strong types extending the given
type $p$ of infinite multiplicity.  The second c.c.c.\ forcing
adds a new automorphism of the algebraic closure of the empty set
that extends to an isomorphism of the models.
In building these models, we 
place a natural measure on the space of strong types extending
$p$ and
introduce a new method of construction.
We require that every element of the construction
realizes a type
over the preceding elements of positive measure.
We expect that this technique can be used to solve other problems within
the context of superstable theories with a type of infinite multiplicity.

In the final section we give a number of examples.  In the first,
we show that $(\RR,\le)$ and $(\RR\ss \{0\},\le)$ are forced isomorphic
by any forcing that adds reals.  
In particular, this shows that the phenomenon of non-isomorphic models
becoming isomorphic in a forcing extension is prevalent, even among
very common structures and forcings as simple as Cohen forcing.
This example also indicates that
membership in a pseudo-elementary class is not absolute, even for
very reasonable forcings.
The second example shows that there is a difference between the
notions of ``potentially isomorphic via c.c.c.'' and ``potentially
isomorphic via Cohen forcing.''  
The third example shows that the assumption of $D(T)$ countable in
Theorem~\ref{big} cannot be replaced by the weaker assumption of $T$
countable.

We assume only that the reader has a basic understanding of stability 
theory and forcing.
On the model theory side, 
all that is required 
is a knowledge of the
basic facts of strong types and forking 
(see \cite{Baldwinbook}, \cite{Lascarbook}, 
\cite{Mak}
or \cite{Shelahbook2nd}).
We assume that our domain of discourse is  a 
large, saturated model $\C$ of $T$.
That is, all models can be taken as elementary submodels of $\C$ and
all sets of elements are subsets of the universe of $\C$.
In Section~2 we work in an expansion $\Ceq$ of $\C$ so that we may consider
strong types to be types over algebraically closed sets.
The definition of $S^+(A,B)$ does not depend on the choice
of the expansion.

Other than a knowledge of the basic techniques of forcing, we 
assume the reader be familiar with the notion of a complete embedding
and basic facts about c.c.c.\ forcings.  The material in 
\cite{Kunen} is more than adequate.

\def\sp{S^*_{p_0,\dots,p_{n-1}}}
\def\Aut{{\rm Aut}}

\section{Strong types and measures}

Throughout this section, assume that $T$ is countable and stable.
As we will be concerned with the space of strong types extending a
given type, it is convenient to fix an expansion $\Ceq$ of $\C$,
where the signature of $\Ceq$ 
contains a sort corresponding to each definable equivalence
relation $E$ of $\C^n$, and a function symbol $f_E$ taking each
tuple to its corresponding $E$-class in its sort.
The advantage of this assumption is that all types are
stationary over algebraically closed sets in $\Ceq$ (see \cite{Mak}).

Our goal in this section is to define a measure on the space of strong
types extending a given type.  Using this measure, we are interested in the
subsets having positive measure.  This leads to our definition
of $S^+(A,B)$.

\begin{defn}  For $p\in S_1(B)$, $B$ finite, let $S^*_p=\{r\in S_1(\acl(B)):p\sbq r\}$.
\end{defn}

As we are working in $\Ceq$,
there is a natural correspondence between 
$S^*_p$ and the set of all strong types extending $p$.
We endow $S^*_p$ with a natural topology $\tau$ by taking as a base
all sets of the form
$$[a/E]=\{r\in S^*_p:r(x)\vdash E(x,a)\}$$
for some equivalence relation $E$ over
$B$ with finitely many classes and some realization $a$ of $p^\C$.

As $T$ is countable and $B$ is finite, there are only countably many equivalence
relations over $B$, so $\tau$ is separable.
In addition $\Aut_B(\C)$ acts naturally on $S^*_p$, so for each
equivalence class $[a/E]$, let ${\rm Stab}([a/E])$ denote the setwise stabilizer
of $[a/E]$.  As $E$ has only finitely many classes, Stab$([a/E])$ has finite index
in $\Aut_B(\C)$.  We construct a regular measure $\mu_p$ on 
the Borel subsets of $S^*_p$ 
by defining $\mu_p([a/E])=1/n$, where $n$ is the index of Stab$([a/E])$ in
$\Aut_B(\C)$ and inductively extending the measure to the Borel subsets.
This is nothing more than the usual construction of Haar measure on
the range of a group action (see e.g., \cite{Haar}).
It is easy to see that the measure $\mu_p$ induces a complete metric on
$S^*_p$, which implies that $S^*_p$ is a Polish space.

For a finite set $A$ and $q\in S_1(A)$, let $\G^q_p=\{r\in S^*_p:q\cup r$ is 
consistent$\}$.  By compactness, $\G^q_p$ 
is a closed, hence measurable subset of $S^*_p$.
For $B\sbq A$ and $A$ finite, let
$$S^+(A,B)=\{q\in S_1(A):q\text{does not fork over} B\text{and} \mu_p(\G^q_p)>0,
\text{where} p=q|B\}.$$

We remark that instead of looking at sets of positive measure,
we could have defined $S^+(A,B)$ to be the set of 
non-forking extensions $q$ of  $p$ such that $\Gamma^q_p$
is non-meagre.  These two notions are not the same,
but they share many of the same properties.  In particular, all of the
lemmas of this section have analogs in the non-meagre context.

\begin{lemma}  \label{product}
Assume $C\sbq B\sbq A$, $A$ finite and that $q\in S_1(A)$ 
does not fork over $C$.  Let $p=q|C$. Then
$\mu_p(\G^q_p)=0$ if and only if either $\mu_p(\Gamma^{q|B}_{p})=0$ 
or $\mu_{q|B}(\Gamma^q_{q|B})=0$.
\end{lemma}

\bp 
For any equivalence class $E$ over $C$ with finitely many classes, say that
$[d/E]$ is consistent with $q$ if there is a realization $e$ of $q$  with
$E(d,e)$.
As $q$ does not fork over $C$, there is a homeomorphism between
$S^*_{q|B}$ and the subspace $\Gamma^{q|B}_p$ of $S_p^*$,
but $\mu_p([d/E])$ may  not equal $\mu_{q|B}([d/E])$.
However, it follows directly from the definitions of the measures that
$\mu_p([d/E])\le\mu_{q|B}([d/E])$ for all $[d/E]$ consistent with $q|B$.
Hence $\mu_p(\Gamma^q_p)\le\mu_{q|B}(\Gamma^q_{q|B})$.
Trivially, $\Gamma_p^q\sbq \Gamma^{q|B}_p$, so
$\mu_p(\Gamma_p^q)\le\mu_p(\Gamma^{q|B}_p)$, which completes the proof of the
lemma from right to left.

For the converse, let $\l=\mu_p(\Gamma^{q|B}_p)$.
We will show that $\l\cdot\mu_{q|B}(\Gamma^q_{q|B})\le\mu_p(\Gamma^q_p)$.
For this, it suffices to show that
$$\l\cdot\mu_{q|B}([d/E])\le\mu_p([d/E])$$
for every $[d/E]$ consistent with $q$.
By definition of the measures, $\mu_p([d/E])=1/n$, where $n$ is the number of
$E$-classes consistent with $p$ and $\mu_{q|B}([d/E])=1/m$, where
$m$ is the number of $E$-classes consistent with $q|B$.
Thus, we must show that $\l\le m/n$.
To see this, let $d_0,\dots,d_{m-1}$ enumerate the $E$-classes consistent with
$q|B$.  As $\bigcup_{i<m}[d_i/E]$ is a disjoint open cover of $\Gamma^{q|B}_p$
and $\mu_p([d_i/E])=1/n$ for each $i$, the regularity of $\mu_p$ implies
that $\l\le m/n$.\nobreak \endproof

\begin{lemma}  \label{transitive}
If $C\sbq B\sbq A$ and $A$ is finite then for every $a$,
$\tp(a/A)\in S^+(A,C)$ if and only if 
$\tp(a,A)\in S^+(A,B)$ and $\tp(a/B)\in S^+(B,C)$.
\end{lemma}

\bp  Let $q=\tp(a/A)$.  As non-forking is transitive,
$q$ does not fork over $C$ if and only if $q$ does not fork over $B$ and
$q|B$ does not fork over $A$.  Further, by Lemma~\ref{product},
$\mu_{q|C}(\G^q_{q|C})>0$ if and only if
$\mu_{q|C}(\G^{q|B}_{q|C})>0$ and 
$\mu_{q|B}(\G^q_{q|B})>0.$
\endproof

Suppose that $p_0,\dots,p_{n-1}\in S_1(B)$.
Let $\sp=\{r\in S_n(\acl(B)):r\r x_i=p_i$ and if $\cbar$ realizes $r$
then $\{c_i:i<n\}$ is independent over $B\}$.

We endow $\sp$ with the analogous topology as $\tau$. 
As types over algebraically closed sets (in $\Ceq$)
have unique non-forking extensions
to any superset of their domain, 
$\sp$ is homeomorphic to the topological
product $\Pi_{i<n} S^*_{p_i}$.
Via this identification, 
endow $\sp$ with the product measure $\mu_{p_0,\dots,p_{n-1}}=
\mu_{p_0}\times\dots\times\mu_{p_{n-1}}$ on the basic open sets and
extend the measure to the Borel subsets.

For $q\in S_n(A)$, let $\G^q_p=\{r\in \sp: r\cup q$ is consistent$\}$.
As before,
$\G^q_p$ is a closed, hence measurable subset of $\sp$.
For $B\sbq A$, $A$ finite, let
$$S^+_n(A,B)=\{q\in S_n(A):q\text{does not fork over} B\text{and}
\mu_{p}(\G^q_p)>0,\text{where} p=q|B\}.$$

The proof of the following lemma is basically an application of Fubini's Lemma
to our context.

\begin{lemma}  \label{Fubini}
Assume that $q(x,y)\in S_2(A)$, $B\sbq A$, $A$ finite and 
that $q$ does not fork over $B$.
Let $q_0=q\r x$, $q_1=q\r y$, and let $p,p_0,p_1$ denote the restrictions of
$q,q_0,q_1$ (respectively) to $B$.
Let $b$
be any realization of $q_1$ 
and let $\G_{p_0}^{q_b}=\{r\in S^*_{p_0}:
r(x)\cup q(x,b)$ is consistent$\}$.
Then
$$\mu_{p_0p_1}(\G_p^q)=\int_{S^*_{p_1}}\mu_{p_0}(\G^{q_b}_{p_0}) d\mu_{p_1}
=\mu_{p_0}(\G^{q_b}_{p_0})\cdot\mu_{p_1}(\G^{q_1}_{p_1}).$$
\end{lemma}

\bp
The first equality is literally Fubini's Lemma and the second follows from
the fact that $\G^{q_b}_{p_0}=\emptyset$ unless $b$ realizes $q_1$ and
the fact that $\mu_{p_0}$ is invariant under translations by elements
of $\Aut_B(\C)$.
\endproof

\begin{lemma} \label{symmetrylemma}
If $B\sbq A$ and $A$ is finite, then for all $a,b$,
$\tp(ab/A)\in S^+_2(A,B)$ if and only if
$\tp(a/A\cup\{b\})\in S^+(A\cup\{b\},B)$ and $\tp(b/A)\in S^+(A,B)$.
\end{lemma}

\bp
This follows from Lemma~\ref{Fubini} in the same manner as Lemma~\ref{transitive}
followed from Lemma~\ref{product}.
\endproof

The following two lemmas are the goals of this section.
The first is the key ingredient in the proof of the Generalized Symmetry Lemma
(Lemma~\ref{GSL}).
The reader should compare it to Axiom VI in [\cite{Shelahbook2nd}, Section IV.1].
The second, the Extendibility Lemma, makes critical use of the added hypothesis
that $|D(T)|=\aleph_0$ that will be assumed throughout the next section.

\begin{lemma}  \label{symmetry}
Assume that $T$ is stable and countable, $B,C\sbq A$, $A$ finite,
$\tp(a/A)\in S^+(A,B)$ and $\tp(b/A\cup\{a\})\in S^+(A\cup\{a\},C)$.
Then $\tp(a/A\cup\{b\})\in S^+(A\cup\{b\},B)$.
\end{lemma}

\bp
Let $D=B\cup C$.  By Lemma \ref{transitive}, $\tp(a/A)\in S^+(A,D)$ and
$\tp(b/A\cup\{a\})\in S^+(A\cup\{a\},D)$.  By Lemma~\ref{symmetrylemma} (switching the
roles of $a$ and $b$),
$\tp(ab/A)\in S^+_2(A,D)$.  Using Lemma~\ref{symmetrylemma} again,
$\tp(a/A\cup\{b\})\in S^+(A\cup\{b\},D)$.  So 
$\tp(a/A\cup\{b\})\in S^+(A\cup\{b\},B)$ using Lemma~\ref{transitive}.
\endproof

\begin{lemma}  [Extendibility Lemma] \label{extendibility}
Assume that $T$ is countable and stable and that $|D(T)|=\aleph_0$.
Let $C\sbq B\sbq A$ be finite, let $E$ be an equivalence relation with finitely
many classes and let $a$ be arbitrary.
If $q\in S^+(B,C)$, $p=q|C$ and $\mu_p([a/E]\cap\Gamma^q_p)>0$,
then there
is $q^+\in S^+(A,C)$ extending $q\cup\{E(x,a)\}$.
\end{lemma}

\bp
Let $\{q_i:i\in\o\}$ enumerate the non-forking extensions of $q$ to $S_1(A)$
that are consistent with $E(x,a)$.
We claim that $[a/E]\cap\G^q_p=\bigcup_{i\in\o}\G^{q_i}_p$.  For, if
$r\in [a/E]\cap\G^q_p$, then as $r\cup q\cup E(x,a)$ is consistent 
we can choose a realization
$b$ of it with $\fg b B A$.  Then $r\in \G^{q_i}_p$, where 
$q_i=\tp(b/A)$.
As $\mu_p$ is countably additive, $\mu_p(\G^{q_i}_p)>0$ for some $i\in\o$.
As non-forking is transitive this $q_i\in S^+(A,C)$, as desired.
\endproof

\section{Positive measure constructions}
\def\var{{\rm var}}

In this section we define two partial orders $(\P,\le_\P)$ and $(\R,\le_\R)$
that will be used in the 
proof of Theorem~\ref{big}.
The forcing $\P$ will force the existence of countable subsets 
$B$ and $C_\al$ ($\al\in\o_1$) of $\C$ such that 
$\acl(B)$ and $\acl(B\cup C_\al)$
are $\ao$-saturated models of $T$, $\acl(B)\preceq\acl(B\cup C_\al)$,
and $\{C_\al:\al\in\o_1\}$ are independent over $B$.
{\it Throughout this section, assume that $T$ is stable,
$|D(T)|=\ao$ (hence $|T|=\ao$)
and we have a fixed type $r^*\in S_1(\emptyset)$ of infinite
multiplicity.}

\begin{defn}  Let $\V=X\cup\bigcup_{\al\in\o_1} Z_\al$, where
$X=\{x_m:m\in\o\}$ and each set $Z_\al=\{z^\al_m:m\in\o\}$,
$\al\in\o_1$
is  a countable set of distinct variable symbols.  
A {\it $\V$-type\/} $q$ is a complete type
in finitely many variables of $\V$.  Let $\var(q)$ denote this set of variables.
\end{defn}

A $\V$-type should be thought of as the type of a finite 
subset of $A\cup\bigcup_{\al\in\o_1} B_\al$.
As notation, given a sequence $\langle a_i:i<n\rangle$ and $u\sbq n$,
let $A_u=\{a_j:j\in u\}$.  Note that as a special case, $A_i=\{a_j:j<i\}$.

\begin{defn}    \label{approximation}
A {\it positive measure construction\/}
(PM-construction) $t$ (of length $n$) is a 
sequence of triples
$\langle (a_i,u_i,v_i):i<n\rangle$ 
satisfying the following conditions for each $i<n$: 
\begin{enumerate}
\item $\tp(a_i/A_i)$ is not algebraic and $u_i\sbq i$;
\item $\tp(a_i/A_i)\in S^+(A_i,A_{u_i})$;
\item If $v_i\in X$ and $j\in u_i$ then $v_j\in X$;
\item If $v_i\in Z_\al$ for some $\al$ and $j\in u_i$ then 
$v_j\in X\cup Z_\al$;
\item If $v_i=z_0^\al$ then $u_i=\emptyset$ and $\tp(a_i/\emptyset)=r^*$.
\end{enumerate}

A PM-construction $t$
may be thought of as
a way of building the $\V$-type
$\tp(a_i:i<n)$ in the
variables $\langle v_i:i<n\rangle$.  
Let $\tp(t)$ denote this type and let $\var(t)=\{v_i:i<n\}$.
If $\tp(t)=q$ then we call $t$ a
{\it PM-construction of $q$.\/}
A $\V$-type $q$ is {\it PM-constructible\/} if there is a PM-construction of it.
\end{defn}

Intuitively, $(a_i,u_i,v_i)\in t$ ensures that
$\tp(a_i/A_i)$ is as generic as possible, given that it
extends $\tp(a_i/A_{u_i})$.
Clause~(5) implies that the set $\{z^\al_0:\al\in\o_1\}\cap\var(t)$
realizes a generic subset of the strong types extending $r^*$.  In 
particular, no two such variables can realize the same strong type.

\begin{defn}  Let $\P$ denote the set of all PM-constructible $\V$-types.
For $p,q\in \P$, say $p\le_{\P} q$ if and only if there is a PM-construction $t$
of $q$ and an $m\in\o$ such that $t\r m$ is a PM-construction of $p$.
That $\le_\P$ induces a partial order on $\P$ follows from the lemma below.
\end{defn}

\begin{lemma}  \label{continue}
Assume that $p\le_{\P} q$.  Then any PM-construction of $p$ can be continued to
a PM-construction of $q$.
\end{lemma}

\bp
Suppose that $t=\langle (a_i,u_i,v_i):i<n\rangle$ is a PM-construction of
$q$ such that $t\r m$ is a PM-construction of $p$ and let
$s=\langle (b_j,u_j',v_j'):j<m\rangle$ be any PM-construction of $p$.
Since $\{v_i:i\in m\}=\{v_j':j\in m\}$ setwise there is a unique permutation
$\sigma$ of $n$ such that $v_i=v_{\sigma(i)}'$ for all $i<m$
and $\sigma(i)=i$ for all $m\le i<n$.
As $\tp(a_i:i<m)=\tp(b_{\sigma(i)}:i<m)$, we can choose an automorphism $\psi$ of
$\C$ such that $\psi(a_i)=b_{\sigma(i)}$ for each $i$.
It is now easy to verify that $s\conc\langle 
(\psi(a_k),\sigma''(u_k),v_k):m\le k<n\rangle$
is a PM-construction of $q$ continuing $s$.
\endproof

The following lemma will be used to show that a generic subset of $\P$
generates a family of  $\ao$-saturated models of $T$.

\begin{lemma}  \label{density}
Let $t=\langle(a_i,u_i,v_i):i<n\rangle$ be any PM-construction.
\begin{enumerate}
\item If $x_m\in X\ss\var(t)$ and $u\sbq n$ such that $j\in u$ implies $v_j\in X$
and $p$ is a non-algebraic 1-type over $A_u$
then there is a realization $a_n$ of $p$ such that
$t\conc\langle(a_n,u,x_m)\rangle$ is a PM-construction.

\item If $z_m^\al\in Z_\al\ss\var(t)$, $m\neq 0$,
$u\sbq n$ such that $j\in u$ implies $v_j\in X\cup Z_\al$ 
and $p$ is a non-algebraic 1-type over $A_u$
then there is a realization  $a_n$ of $p$ such that
$t\conc\langle(a_n,u,z_m^\al)\rangle$ is a PM-construction.

\item If $z_0^\al\in Z_\al\ss\var(t)$,
then there is an $a_n$ such that
$t\conc\langle(a_n,\emptyset,z_0^\al)\rangle$ is a PM-construction.
\end{enumerate}
\end{lemma}

\bp
These follow immediately from the Extendibility Lemma
and Clauses~(3),~(4),~(5) of Definition~\ref{approximation}.
\endproof

In order to establish the independence of the $B_\alpha$'s over $A$
and to analyze the complexity of the partial order $(\P,\le_\P)$, 
we seek a `standard form' for a PM-construction.
The primary complication 
is that the restriction of a PM-constructible type to a subset of 
its free variables need not
be PM-constructible.  
We characterize when a permutation $\sigma$
of a PM-construction $t$ is again a PM-construction.  
Call a permutation $\sigma$ {\it permissible\/} if 
$\sigma''(u_i)\sbq \sigma(i)$ for all $i<n$.
Clearly, if $\sigma$ is not permissible then $\sigma t$ violates Clause~(1)
of being a PM-construction.
The following lemma, known as the Generalized Symmetry Lemma, establishes the converse.
Its proof simply amounts to bookkeeping once we have Lemma~\ref{symmetry}.

\begin{lemma}[Generalized Symmetry Lemma]  \label{GSL}
If $t$ is a PM-construction of $q$  and $\sigma$ is a permissible
permutation then $\sigma t$ is a PM-construction of $q$ as well.
\end{lemma}

\bp  Suppose that $t=
\langle (a_i,u_i,v_i):i<n\rangle$ is a PM-construction of $q$.  Then
Lemma~\ref{symmetry} insures that
$\sigma_k(t)$ is a 
PM-construction, where $\sigma_k$ is the (permissible)
permutation exchanging $k$ and $k+1$ whenever
$k\neq n-1$ and $k\not\in u_{k+1}$.
The lemma now follows easily by induction on the length of $t$.
The reader is encouraged to compare this with \cite[IV, Theorem 3.3]{Shelahbook2nd}.
\endproof

As an application of Lemma~\ref{GSL}, 
we obtain a `standard form' for a PM-construction.
Given any $p\in\P$, there is a PM-construction
$t=\langle(a_i,u_i,v_i):i<n\rangle$ of $p$ such that, for all $i<j<n$,
\begin{enumerate}
\item if $v_j\in X$ then $v_i\in X$;
\item if $v_i\in Z_\al$ and $v_j\in Z_{\al'}$ then $\al\le\al'$; 
\item if $v_j= z_0^\al$ then $v_i\not\in Z_\al$.
\end{enumerate}
To see this, let $s$ be any PM-construction of $p$, find an appropriate permissible
permutation $\sigma$ and let $t=\sigma s$.
The following lemma is a consequence of this representation.

\begin{lemma}  \label{nonforking}
Let $p(\xbar,\zbar^{\al_0},\dots,\zbar^{\al_{k-1}})\in\P$, where $\xbar\sbq X$
and $\zbar^{\al_i}\sbq Z_{\al_i}$ and let $\bbar\cbar^{\al_{0}}\dots\cbar^{\al_{k-1}}$
realize $p$.  Then $\{\cbar^{\al_i}:i<k\}$ is independent over $\bbar$.
\end{lemma}

\bp
We argue by induction on $\var(p)$.  Choose $p\in\P$ with $n+1$ free variables.
We can find a PM-construction 
$t=\langle(a_i,u_i,v_i):i<n+1\rangle$ of $p$ with the variables arranged as in the
application above.  
By elementarity, we may assume that $\abar=\bbar\cbar^{\al_{0}}\dots\cbar^{\al_{k-1}}$.
If $v_n\in X$ there is nothing to prove, so say $v_n\in Z_{\al_{k-1}}$
and let $\dbar=\cbar^{\al_{k-1}}\ss \{a_n\}$.
\relax From our inductive hypothesis, 
$\{\cbar^{\al_i}:i<k-1\}\cup\{\dbar\}$ is independent over $\bbar$.
In particular, ${\fg \dbar \bbar
{\{\cbar^{\al_i}:i<k-1\}}}$.  However, $\tp(a_n/A_n)\in S^+(A_n,A_{u_n})$ and
$A_{u_n}\sbq \bbar\cup\dbar$,
so $\tp(a_n/A_n)$  does not fork
over $\bbar\cup\dbar$.  Hence,
$\{\cbar^{\al_i}:i<k\}$ is independent over $\bbar$ by the transitivity of non-forking.
\endproof


\begin{lemma}  \label{upper}
Assume that $p,q_1,q_2\in\P$, $p\le_\P q_1$, $p\le_\P q_2$
and $var(q_1)\cap\var(q_2)=\var(p)$.  Then there is an upper bound $p^*\in\P$ 
of both $q_1$ and $q_2$.
\end{lemma}

\bp
Say $|\var(p)|=n_0$.
Let $s$ be any PM-construction for $p$ and, using Lemma~\ref{continue},
let $t_1=\langle(a_i,u_i,v_i):i<n_1\rangle$ and $t_2=\langle
(b_i,u_i',v_i'):i<n_2\rangle$ be PM-constructions for
$q_1,q_2$ respectively, each continuing $s$.
We form a PM-construction $t^*$ by concatenating a `copy' of $t_2\ss s$ to $t_1$.
More formally, let $\dbar=\langle a_i:i<n_0\rangle$ and for each $k$, $n_0\le k<n_2$,
let $u_k''=(u_k'\cap n_0)\cup\{j+(n_1-n_0):j\in u_k'\cap(n_2\ss n_0)\}$.
Using the Extendibility Lemma, we can successively
find a sequence $\langle c_k:n_0\le k<n_2\rangle$
such that $t^*=t_1\conc\langle(c_k,u_k'',v_k'):n_0\le k<n_2\rangle$ is a PM-construction
and $\tp(\dbar\cbar)=q_2$.
Let $p^*$ be the $\V$-type generated by $t^*$.  Visibly, $q_1\le_\P p^*$.
That $q_2\le_\P p^*$ follows from Lemma~\ref{GSL} by taking the permissible permutation
of $t^*$ exchanging $t_1\ss s$ and the copy of $t_2\ss s$.
\endproof

By using the full strength of the Extendibility Lemma,
using the notation in the proof above, if $E$ is an equivalence relation with
finitely many classes, we may further require that
$E(v_i,v_j')\in p^*$ if and only if $\C\models E(a_i,b_j)$.  
This improvement will be crucial in
the proof of Claim~3 of Lemma~\ref{noiso}.

A partially ordered set $\P$ has the {\it Knaster condition\/} if,
given any uncountable subset $X$ of $\P$, one can find an uncountable $Y\sbq X$ such
that any two elements of $Y$ are compatible.
Evidently, if a partially ordered set has the Knaster condition, then it satisfies
the countable chain condition (c.c.c.).  However, in contrast to the
case for c.c.c.\ posets, it is routine to check that
the product of two posets with the Knaster condition must have the Knaster condition.

\begin{lemma}  \label{Knaster}
$(\P,\le_\P)$ satisfies the Knaster condition, hence $\P\times\P$ satisfies the
countable chain condition.
\end{lemma}

We begin with a combinatorial lemma that is of independent interest.
It is not claimed to be new, but the authors know of no published reference.

\begin{lemma}
There is a partition of 
$[\o_1]^{<\o}$ into $\bigcup\{A_i:i\in\o\}$ such that $c\cap d$
is an initial segment of $c$ whenever
$i\in\o$ and $c,d\in A_i$.
\end{lemma}

\bp
Clearly it suffices to partition each $[\o_1]^l$, so fix $l\in\o$.
We define three families of functions.  First, for each
$\b\in\o_1$, choose an injective function $g_\b:\b\ra\o$.  
Next, for each $n\in\o$, define a partial function $f_n:\o_1\ra\o_1$
by $f_n(\b)=g_\b^{-1}(n)$.  Note that if $f_n(\b)$ is defined, then it is less than
$\b$.  So define a (total) function $h_n:\o_1\ra\o$, where
$h_n(\b)$ is the least $m$ such that the $m$-fold composition $f^{(m)}_n(\b)$ is
undefined.

Define an equivalence relation $\sim$ on $[\o_1]^l$ by putting
$\{\al_0,\dots,\al_{l-1}\}\sim\{\b_0,\dots,\b_{l-1}\}$
if and only if, for all $i<j<l$ there is an $n\in\o$ such that
$n=g_{\al_j}(\al_i)=g_{\b_j}(\b_i)$ and $h_n(\al_j)=h_n(\b_j)$. 
Clearly, $\sim$ partitions $[\o_1]^l$ into countably many classes, so let
each $A_i$ denote a $\sim$-class.

To see that this works, suppose 
$c=\{\al_0,\dots,\al_{l-1}\}\sim d=\{\b_0\dots,\b_{l-1}\}$.
We first observe that if $\al_i=\b_j$ then $i=j$.
For, if not, we could assume by symmetry that $i<j$.  
Let $n=g_{\al_j}(\al_i)=g_{\b_j}(\b_i)$.  Now, $f_n(\al_j)=\al_i=\b_j$,
hence $h_n(\al_j)=h_n(\b_j)+1$, contradicting $c\sim d$.
Next, suppose $\al_j=\b_j$ for some $j<l$ and fix $i<j$.
As $c\sim d$, 
$g_{\al_j}(\al_i)=g_{\b_j}(\b_i)$, so $\al_i=\b_i$ as $g_{\al_j}$ is injective.
Hence $c\cap d$ is an initial segment of both $c$ and $d$.
\endproof

\par{\bf Proof of Lemma~\ref{Knaster}.}$\ \ $
As notation, for each $p\in\P$,
let $u^p$ denote the (finite) set of all $\al$ such that
$\var(p)\cap Z_\al\neq\e$.
Given $f$ a permutation of $\o_1$, $f$ induces a permutation of $\V$
(also called $f$), defined by $f(x_m)=x_m$ and $f(z^\al_m)=z^{f(\al)}_m$.
This $f$ induces a permutation of $\P$, where
$$f(p)=
\{\phi(f(v_0),\dots,f(v_{n-1})):\phi(v_0,\dots,v_{n-1})\in p\}.$$

\relax From the lemma above and the fact that $D(T)$ is countable,
it is easy to find a partition 
$\P=\bigcup\{P_n:n\in\o\}$ such that, for each $n\in\o$ 
and each $p,q\in P_n$,
\begin{enumerate}
\item $|u^p|=|u^q|$;
\item $u^p\cap u^q$ is an initial segment of both $p$ and $q$;
\item if $f$ is any permutation of $\o_1$ fixing $u^p\cap u^q$ pointwise and
$f''(u^p)=u^q$, then $f(p)=q$.
\end{enumerate}

We claim that every pair $p,q\in P_n$ are compatible.  For, let $w=u^p\cap u^q$ and
let $\V_w=X\cup\bigcup\{ Z_\al:\al\in w\}$.  By Clause~(3), 
$\var(p)\cap \V_w=\var(q)\cap\V_w$.  Further, letting $p_0=p\r (\var(p)\cap\V_w)$,
it follows from Clause~(3) and the standard form following Lemma~\ref{GSL}
that $p_0\le_\P p$ and $p_0\le_\P q$.
Thus, by Lemma~\ref{upper}, $p$ and $q$ are compatible.
\endproof

We next define our second forcing notion, $\R$.  We begin by defining
a dense suborder of $\R$.  The intuition behind a faithful triple
$(p,q,h)$ is that $p$ and $q$ are finite approximations to  models of $T$ 
and $h$
is an elementary map between the approximations.

\begin{defn}
A triple $(p,q,h)$ is {\it faithful\/} if $p,q\in\P$ and $h:\var(p)\ra\var(q)$
satisfy:
\begin{enumerate}
\item $h$ is onto;
\item $\phi(\vbar)\in p$ if and only if $\phi(h(\vbar))\in q$ for all 
formulas $\phi(\xbar)$;
\item for $v\in\var(p)$, $v\in X$ if and only if $h(v)\in X$;
\item for $v\in\var(p)$, $\al\in\o_1$, $v\in Z_\al$ if and only if $h(v)\in Z_\al$.
\end{enumerate}
\end{defn}

\begin{lemma}  \label{faithextend}
Suppose $(p,q,h)$ is faithful and $p\le_\P p'$.  There is $q'\ge_\P q$ and
$h'\contains h$ such that $(p',q',h')$ is faithful.
Further, for any finite $F\sbq \V$, we may assume $\var(q')\cap F\sbq\var(q)$.
\end{lemma}

\bp
Let $m=|\var(p)|$.  Arguing by induction on the size of the difference,
we may assume that $|\var(p')|=m+1$.
Let $s'=\langle (a_i,u_i,v_i):i\le m\rangle$ be a PM-construction of $p'$ such
that $s= s'\r m$ is a PM-construction of $p$ and let
$t=\langle (b_j,u_j,w_j):j<m\rangle$
be a PM-construction of $q$.  Our  $h$ induces a map $h^*:A_m\ra B_m$
by putting $h^*(a_i)=b_j$, where $h(v_i)=w_j$.  As 
$(p,q,h)$ is faithful, $h^*$ is elementary.
Let $p_m=\tp(a_m/A_m)$ and let $q_m=h^*(p_m)$.
By elementarity, $q_m$ is a non-algebraic 1-type over $B_m$.
Pick $w_m\in\V\ss(\var(q)\cup F)$ 
such that $w_m\in X$ if $v_m\in X$ and $w_m\in Z_\al$, where
$v_m\in Z_\al$, otherwise.
By Lemma~\ref{density}, there are $b_m$ and $u_m$ such that 
$t'=t\conc\langle(b_m,u_m,w_m)\rangle$
is a PM-construction and $\tp(b_m/B_m)=q_m$.  Thus, $(p',q',h')$ is faithful,
where $q'=\tp(t')$ and $h'=h\cup\{(v_m,w_m)\}$.
\endproof

\begin{defn}
$\R=\{(p,q,h):$  there are $p_1\le_\P p$ and $q_1\le_\P q$ such that $(p_1,q_1,h)$ is
faithful$\}$.
Define a preorder $\le_0$ on $\R$ by $(p,q,h)\le_0 (p',q',h')$ if and only if
either $p\le p'$, $q\le q'$ and 
$h=h'$; or $p=p'$, $q=q'$  $h\sbq h'$ and $(p,q,h')$ is faithful.
Let $\le_\R$ be the transitive closure of $\le_\R$.  It is clear that $\le_\R$ is
a partial order on $\R$.
\end{defn}

\begin{lemma}  \label{faithdense}
The set of faithful triples is a dense suborder of $\R$.
\end{lemma}

\bp
Pick $(p,q,h)\in\R$ and assume that $h:\var(p_1)\ra \var(q_1)$, where $p_1\le_\P p$
and $q_1\le_\P q$.  
By Lemma~\ref{faithextend} there is $q_2\ge_\P q_1$ with $\var(q_2)\cap\var(q)=\var(q_1)$
and $h_2\contains h$ such that $(p,q_2,h_2)$ is faithful.
By Lemma~\ref{upper} there is an upper bound $q^*\in \P$ of both $q_2$ and $q$.
Now consider the triple $(p,q^*,h_2)$.
By Lemma~\ref{faithextend} again (with the roles of $p$ and $q$ reversed)
there is $p^*\ge_\P p$ and $h^*\contains h_2$ such that
$(p^*,q^*,h^*)$ is faithful.  Also, $(p,q,h)\le_0 (p^*,q^*,h)\le_0 (p^*,q^*,h^*)$,
so $(p,q,h)\le_\R (p^*,q^*,h^*)$.
\endproof

\begin{lemma}  \label{natural}
The natural embedding $i:\P\times\P\ra \R$, defined by
$i(p,q)=(p,q,\e)$ is a complete embedding.
\end{lemma}
 
\bp
Fix  a maximal antichain $A\sbq \P\times \P$.
We must show that $i''(A)$ is a maximal antichain in $\R$.
So, fix $(p,q,h)\in\R$.
Choose an element $(p',q')\in\P\times\P$ that is an upper bound of
$(p,q)$ and some $(p_0,q_0)\in A$.
By Lemma~\ref{faithdense}, there is a faithful triple $(p^*,q^*,h^*)\ge_\R (p',q',h)$.
It is easy to check that $(p^*,q^*,h^*)$ is an upper bound of both $(p,q,h)$ and
$i(p_0,q_0)$.
\endproof
\section{The main theorem}

This section is devoted to proving the following theorem.

\begin{theorem}  \label{big}
Assume that $T$ is superstable, $|D(T)|=\ao$ and there is a type of infinite
multiplicity.
There is a c.c.c.\ partial order $\Q$ such that \break
$\forces_\Q$ ``There are two non-isomorphic, potentially isomorphic models of $T$.''
\end{theorem}

\begin{remark}
The forcing $\Q$ will be $\P\times\P$ from the last section, which, in addition to
having the c.c.c., satisfies the Knaster condition.
The second forcing (i.e., $\R/H$)
is almost an element of the ground model $V$.
That is, the forcing $\R\in V$ and
for any $\Q$-generic filter $H$, $\R/H$ will be a  partial order 
forcing the two models isomorphic.
\end{remark}

\bp
By Theorem 1.13 of \cite{BaldLasSheforceiso}, 
we may assume that in addition, $T$ has NDOP
and NOTOP.  In particular, prime and minimal models exist over independent
trees of models of $T$.
Given an $n$-type of infinite multiplicity,
one can find a finite set $\abar$ and a 1-type 
$r^*\in S_1(\abar)$
of infinite multiplicity.  
Let $T'$ denote the $L(\abar)$-theory of $(\C,\abar)$.
Thus, working with $T'$ as our basic theory, $r^*$ is a type over the empty set,
so our results from Section~3 apply.

Fix $\P$ and $\R$ from Section~3 and let $\Q=\P\times\P$.
By Lemma~\ref{Knaster} $\Q$  satisfies 
the c.c.c., and by Lemma~\ref{natural}
the natural embedding of $\Q$ into $\R$ is a complete embedding.
If $H$ is $\Q$-generic, then $R/H$ embeds naturally into the set of
finite partial functions $f:\o_1\times\o\ra\o_1\times\o$ that fix
the first coordinate.  Thus $R/H$ satisfies the c.c.c.
In the remainder of the section we show that $\Q$ `constructs' two non-isomorphic
models and that $\R/H$ forces them isomorphic.

We first show that the forcing $\P$ constructs a new model of our theory,
i.e., one that is not isomorphic to any structure in the original universe $V$.
Fix a $\P$-generic filter $G$.
We associate a model $\B^*[G]$ of $T$ with $G$ as follows.
First, by applying Lemma~\ref{density}, for every $\vbar\sbq \V$, $\{p\in\P:
\vbar\sbq \var(p)\}$ is dense, hence there is a $p\in G$ such that $\vbar\sbq \var(p)$.
In addition, as any $p,q\in G$ have a common upper bound, $p\r\vbar=q\r\vbar$
for any $\vbar\sbq\var(p)\cap\var(q)$.
Let$$\Gamma_G=\{\phi(\vbar):\vbar\sbq\V\ \hbox{and}\ \phi(\vbar)\in p\ \hbox{for some}
\ p\in G\}.$$
Let $A=A_X\cup\bigcup_{\al\in\o_1}A_\al$ be a realization of
$\Gamma_G$ in $\C$  i.e., for all
$\bbar\sbq A$, $\C\models\phi(\bbar)$ if and only if $\phi(\vbar)\in \Gamma_G$,
where $\vbar$ is the tuple from $\V$ corresponding to $\bbar$.
Working inside $\C$, let 
$\A_\e[G]=\acl(A_X)$ and for each $\al\in\o_1$, $\A_\al[G]=\acl(A_X\cup A_\al)$.

We first claim that $\A_\e[G]$ and each $\A_\al[G]$ is an $\ao$-saturated model of $T$.
To see that this holds of $\A_\e[G]$, 
note that by Lemma~\ref{density}(1), $A_X$ realizes every non-algebraic 1-type over
a finite subset of itself.  It is a straightforward exercise to show that this fact, 
together with $\A_\e[G]=\acl(A_X)$ implies that $\A_\e[G]$ is an $\ao$-saturated
model of $T$.
The proof for each $\A_\al[G]$ is analogous, using Lemma~\ref{density}(2)
to show that $A_X\cup A_\al$ realizes every non-algebraic 1-type over a finite
subset of itself.

Also, it follows from Lemma~\ref{nonforking} that $\{\A_\al[G]:\al\in\o_1\}$ 
is independent over $\A_\e[G]$.
As $T$ satisfies NDOP and NOTOP, we can form a continuous, increasing chain
of countable models 
$\langle \B_\al[G]:\al\in\o_1\rangle$ such that $\B_\al[G]$ is prime and minimal
over $\bigcup_{\b<\al} \A_\b[G]$.
By replacing the chain by an isomorphic copy, we may assume that the universe $B_\al$
of each $\B_\al[G]$ is a countable subset of $\o_1$ and that 
$\B^*[G]=\bigcup_{\al\in\o_1} \B_\al[G]$ has universe $\o_1$.
The crucial fact is that this model $\B^*[G]$ is not $L$-isomorphic to 
any structure in the ground universe.

\begin{lemma}  \label{noiso}
In $V[G]$ there is no $L$-elementary embedding of $\B^*[G]$ into any model $\D\in V$.
\end{lemma}

\bp
Fix $\D\in V$.  As $\C$ is sufficiently saturated, we may assume that $\D$ is an
elementary substructure of $\C$.
Assume by way of contradiction that
such an embedding $f$ exists.
\relax From our assumption above, $f$ is an elementary map between two subsets of $\C$.
Fix a $\P$-name $\f{f}$ and a condition $g_0\in G$
such that 
$$g_0\forces \f{f}:\f{\B^*[G]}\ra D.$$
Also, as $\{\B_\al:\al<\o_1\}$
is a continuous, increasing chain we can find $\P$-names $\f{B_\al}$
such that $\al<\b$ implies $\f{B_\al}\sbq \f{B_\b}$ and 
$\f{B_\d}=\bigcup\{\f{B_\al}:\al<\d\}$.
Further, for any $\al<\o_1$, since $\P$ satisfies c.c.c.\ and
$\forces \hbox{$\f{B_\al}$ is countable}$, there is $\b<\o_1$ such that
$\forces \f{B_\al}\sbq \b.$
Consequently, we may assume that each $\f{B_\al}$ is a countable $\P$-name.

For each $\d\in\o_1$, let 
$\V_\d=X\cup\bigcup_{\al<\d} Z_\al$, let 
$\P_\d=\{p\in\P:\var(p)\sbq \V_\d\}$,
and let $G_\d=G\cap P_\d$.
\Claim{Claim 2}  For all $\d\in\o_1$, the identity map $i:\P_\d\ra\P$ is a complete
embedding.
\smallskip
\bp  
Let $A$ be a maximal antichain in $\P_\d$ and let $p\in \P$.
Let $p_0=p\r\V_\d$.  As $A$ is maximal, there is $q_0\in A$ and $q\in \P_\d$,
$q$ an upper bound of both $p_0$ and $q_0$.
By Lemma~\ref{upper}, there is an upper bound of both $p$ and $q$ (hence of $q_0$).
Thus, $A$ is a maximal antichain of $\P$ as well.
\endproof

Let $\P/G_\d=\{p\in \P:$ $p$ is compatible with each $g\in G_\d\}$
and let $G_\d^*$ be the $\P/G_\d$-generic filter induced by $G$.
As the identity is a complete embedding, 
$V[G]=V[G_\d][G_\d^*]$ (see e.g., \cite{Kunen}).
It is easily verified that a condition $p\in\P$ is an element of $\P/G_\d$ if and
only if $p\r\V_\d\sbq \Gamma_{G_\d}$.
Let $$C=\{\d<\o_1:\f{B_\d}\ \hbox{is a $\P_\d$-name}, \f{B_\d}\sbq\d,
\ \hbox{for all}\ \al<\d, \f{f(\al)}
\ \hbox{is a $\P_\d$-name}\}.$$
Visibly, $C\in V$.  Using the fact that $\P$ satisfies c.c.c.\ again,
$C$ is a club subset of $\o_1$.
Note that $\B_\d[G]\in V[G_\d]$ and $f\r \d\in V[G_\d]$ for each $\d\in C$.
Fix an element $\d\in C$.

\Claim{Claim 3}  Working in $V[G_\d]$,
for each $e\in \D$, the set
$$D^*_e=\{p^*\in \P/G_\d: p^*\forces_{\P/G_\d} 
\tp(z_0^\d,\B_\d[G])\neq \tp(e,f(\B_\d[G])\big)\}$$
is dense in $\P/G_\d$.
\smallskip
\bp
Fix $e$ and choose $p\in \P/G_\d$.
By Lemma~\ref{density} and our characterization of $\P/G_\d$
we may assume $z_0^\d\in\var(p)$.
Let $p_0=p\r\P_\d$ and let $m=|\var(p_0)|$.
Let $t=\langle (a_i,u_i,v_i):i<n\rangle$ be a PM-construction of $p$
such that $t_1=_{{\rm def}}t\r m$ is a PM-construction of $p_0$ and $v_m=z_0^\d$.
Let $\psi$ be an automorphism of $\C$ fixing $A_m$ pointwise such that
$\stp(a_m)\neq\stp(\psi(a_m))$.  
(One exists since $r^*=\tp(a_m/\e)$ has infinite multiplicity and 
$\tp(a_m/A_m)\in S^+(A_m,\e)$.)
Fix a definable equivalence relation
$E$ with finitely many classes such that 
$\C\models\neg E(a_m,\psi(a_m))$ and pick a set of 
representatives $\{c_i:i<k\}$ of $E$'s classes from $\B_\d[G]$.
Say $E(e,f(c_i))$ holds in $\C$.
Choose $g\in G_\d$ such that $g \forces_{\P_\d} E(e,f(c_i))$ and $\var(p_0)\cup\{c_i\}
\sbq \var(g)$.
Let $s$ be a PM-construction of $g$ and suppose $(b,u,c_i)\in s$.
We may assume that $\C\models\neg E(a_m,b)$, since otherwise we could replace
$a_m$ by $\psi(a_m)$ in the argument below.
As in the proof of Lemma~\ref{upper} (and the remark following the proof) 
it 
follows from the Extendibility Lemma that there is a sequence  $s^*$  and an 
element $d$ such that
$s\conc s^*$ is a PM-construction,   
$\tp(s^*)=\tp(t\ss t_1)$,
$(d,\e,z^\d_0)\in s^*$, and $\C\models E(a_m,d)$.
Let $p^*=\tp(s\conc s^*)$.  
As $E$ is an equivalence relation,
$\C\models \neg E(d,b)$ so $\neg E(z_0^\d,c_i)\in p^*$.
Further,
$p^*\r \V_\d=g$, $p^*\in\P/G_\d$ as required.
\endproof

Thus, working in  $V[G]=V[G_\d][G^*_\d]$,
Claim~3 implies that 
$\tp(z_0^\d,\B_\d[G])\neq \tp(e,f(\B_\d[G]))$ for all $e\in D$, contradicting
the elementarity of $f$.
\endproof

Continuing with the proof of Theorem~\ref{big},
fix $H=G_1\times G_2$, a  $\P\times \P$-generic filter.  Following the
procedure above, we can build
elementary substructures $\B^*[G_1]$ and $\B^*[G_2]$ of $\C$ in $V[H]$.
It follows from Lemma~\ref{noiso} and the fact that $V[H]=V[G_1][G_2]$
that there is no $L$-isomorphism 
$f:\B^*[G_1]\ra\B^*[G_2]$ in 
$V[H]$.

To complete the proof of the
theorem, it remains to show that $\B^*[G_1]$ can be forced 
isomorphic to $\B^*[G_2]$ by a c.c.c. forcing.
Let $\R/H=\{(p,q,h)\in\R: (p,q,h)$ is compatible with $i(p',q')$ for
every $(p',q')\in H\}$.
As noted above, $\R/H$ satisfies the c.c.c.\
We claim that $\R/H$ forces  an $L(\abar)$
isomorphism between $\B^*[G_1]$ and $\B^*[G_2]$.
Indeed, let
$$h^*=\bigcup\{h:(p,q,h)\in\R/G_1\times G_2\ \hbox{for some}\ p,q\in\P\}.$$

By Lemma~\ref{faithextend}, 
$h^*$ is an $L(\abar)$-elementary map from a set of realizations
of $\Gamma_{G_1}$ to a set of realizations of $\Gamma_{G_2}$.
Now $h^*$ easily extends to an $L(\abar)$-elementary map of the algebraic closures
of these sets.  That is, $h^*$ maps the independent tree 
$\bigcup\{\A_\al[G_1]:\al\in\o_1\}$ 
of models of $T'$ to  the independent tree
$\bigcup\{\A_\al[G_2]:\al\in\o_1\}$.
As the prime and minimal model of such a tree is unique, $h^*$ extends to an
$L(\abar)$-isomorphism of $\B^*[G_1]$ and $\B^*[G_2]$.
\endproof

\section{Some examples}
Our first example demonstrates the ubiquity of the phenomenon of non-isomorphic
models becoming isomorphic in a forcing extension.
It implies that even very weak forcings such as Cohen forcing
are able to alter the isomorphism type of some very simple structures.

\begin{example}  \label{reals}
Let $M_1=(\RR^V,\le)$ and $M_2=(\RR^V\ss\{0\},\le)$.
Then $M_1$ is not isomorphic to $M_2$ in the ground universe
$V$, but $M_1$ and $M_2$ become 
isomorphic in any transitive $V'\contains V$ with $\RR^{V'}\neq\RR^V$.
\end{example}  

\bp
It is clear that $M_1$ and $M_2$ are not isomorphic in $V$.
Fix $V'$, a transitive extension of $V$ that adds reals.
We will construct an isomorphism $f\in V'$ between $M_1$ and $M_2$.
Towards this end, first note that as $V$ and $V'$ are both transitive,
$\o$, $\ZZ$ and $\QQ$ are all absolute between $V$ and $V'$.
In particular, $\QQ^{V'}=\QQ^V$.  As $\RR^{V'}$ is defined as the set of all
Dedekind cuts of rationals, $\QQ^V$ is dense in $\RR^{V'}$.
   
Next, for any $a,b\in\RR^{V'}$ with $a<b$, fix $\{x_n:n\in\ZZ\}$, a strictly
increasing sequence from $(a,b)$ that is both cofinal and coinitial in $(a,b)$.
Using the density of $\QQ^V$ in $\RR^{V'}$, we may successively
choose $y_n\in\QQ^V\cap(x_n,x_{n+1})$ to 
obtain a cofinal, coinitial sequence  of order-type $\ZZ$
in $(a,b)$ with each element in $V$.

Using, this, we claim that if $a<b$ and $c<d$, then there is an isomorphism
$g:(a,b)\cap\RR^V\ra (c,d)\cap\RR^V$.
To see this, choose strictly increasing sequences $\langle y_n:n\in\ZZ\rangle$
and $\langle z_n:n\in\ZZ\rangle$
from $\QQ^V$, cofinal and coinitial in $(a,b)$ and $(c,d)$, respectively.
Now, as $(y_n,y_{n+1})\cap \RR^V$  and $(z_n,z_{n+1})\cap \RR^V$
are each  open intervals in $\RR^V$, there is an isomorphism $g_n\in V$
between them.  Piecing these isomorphisms
together yields an isomorphism between $(a,b)\cap\RR^V$
and $(c,d)\cap\RR^V$.

We are now ready to build our isomorphism between $M_1$ and $M_2$.
Fix $a<b<c$ in $\RR^{V'}\ss\RR^V$ with $a<0<c$.
\relax From the paragraph above, let $g_1$ be an isomorphism between $(a,b)\cap\RR^V$ and 
$(a,0)\cap\RR^V$ and let 
$g_2$ be an isomorphism between $(b,c)\cap\RR^V$ and 
$(0,c)\cap\RR^V$.
Define $f:M_1\ra M_2$ by
$$f(x)=\left\{ \begin{array}{ll}
                    x & \mbox{if $x<a$ or $c<x$} \\
                    g_1(x) & \mbox{if $a<x<b$} \\
                    g_2(x) & \mbox{if $b<x<c$.}
               \end{array} \right.$$

\par\bigskip

The (pseudo-elementary) class $\Khom$ of homogeneous linear orders is the class of
all dense linear orders with no endpoints such that any non-empty open interval
is isomorphic to the entire linear order.
Examples include $(\QQ,\le)$ and $(\RR,\le)$.
It is well known that forcing preserves satisfaction for models.  Thus, the
relation ``$M\in {\bf K}$'' is absolute between transitive models of set theory
for elementary classes {\bf K}.  Similarly, if {\bf K} is a pseudo-elementary
class (i.e., a class of reducts of an elementary class)
and $M\in {\bf K}$ in the ground universe, then $M\in {\bf K}$ in any forcing
extension.  However, Example~\ref{reals} indicates that the converse need not hold.
That is, $M_2\not\in\Khom$ in $V$, while $M_2\in\Khom$ in any transitive
$V'\contains V$ that adds reals.

The class $\Khom$ can also be used to show that `potential isomorphism
via c.c.c.\ forcings' is distinct from `potential isomorphism via Cohen forcings.'
As $\Khom$ is unstable, it follows from Theorem~1.7 of
\cite{BaldLasSheforceiso} that there is a pair
of non-isomorphic structures in $\Khom$ that can be forced isomorphic by a c.c.c.\
forcing.  This contrasts with the theorem below.

\begin{theorem}
Let $\Q=(^{<\o}\o,\triangleleft)$ be Cohen forcing.  
For all $M_1,M_2\in\Khom$,\break
$M_1\cong M_2$ if and only if
$\forces_\Q M_1\cong M_2$.
\end{theorem}

\bp
Right to left is clear by absoluteness.
Choose homogeneous linear orders $M_1=(I_1,\le)$ and $M_2=(I_2,\le)$ such that
$\forces_\Q M_1\cong M_2$.  We will construct an isomorphism $g:M_1\ra M_2$
in the ground universe as a countable union of approximations in the sense
described below.

Fix a $\Q$-name $\ff$ such that $\forces_\Q$ ``$\ff$ is an isomorphism between $M_1$
and $M_2$.''  For each $q\in\Q$, let $f_q=\{(a,b)\in I_1\times I_2: q\forces \ff(a)=b\}$.
To ease notation, let $I_i'=I_i\cup\{-\infty,\infty\}$ ($i=1,2$), where
$-\infty$ is the smallest element of $I_i'$ and $\infty$ is the largest.
For $h$ a partial 1-1 function from $I_1'$ to $I_2'$, let $D_1(h)=\dom(h)$
and $D_2(h)=\dom(h^{-1})$.

An {\it approximation on $[x_0,x_1]$\/} 
is a partial, order-preserving function $g:[x_0,x_1]\ra I_2'$
such that, for each $a\in [x_0,x_1]\ss D_i(g)$ there are $b,c\in D_i(g)$ with $b<a<c$
and $(b,c)\cap \D_i(g)=\e$.
If $[x_0,x_1]=I_1'$, $g$ is simply called an approximation.

Trivially, $g_0=\{(-\infty,-\infty),(\infty,\infty)\}$ is an approximation.
As noted above, we will construct an increasing sequence $\langle g_n:n\in\o\rangle$
of approximations such that for each $q\in\Q$ there is $n\in\o$ such that
$D_i(g_n)\contains D_i(f_q)$ ($i=1,2$).  Once we build
such a sequence, $g=\bigcup g_n$ will be an isomorphism between $I_1'$
and $I_2'$ since every $a\in I_1$ is in $D_i(f_q)$ for some $q\in \Q$.
Thus, all that remains is to prove the following claim.

\Claim{Claim}
For every approximation $g$ and $q\in \Q$ there is an approximation $g'\contains q$
with $D_i(g')\contains D_i(f_q)$, $i=1,2$.

\bp
Fix an approximation $g$ and $q\in\Q$.  By symmetry it suffices to find $g'\contains g$
with $D_1(g')\contains D_1(f_q)$.
As $D_1(g)$ partitions $I_1$ into convex sets, we may independently 
find approximations $g'$ on $[x_0,x_1]$ extending $g\r[x_0,x_1]$
for each pair $x_0,x_1\in D_1(g)$ with 
$x_0<x_1$ and $(x_0,x_1)\cap D_1(g)=\e$.
So fix such a pair $(x_0,x_1)$.  Choose $p\ge q$ such that $x_0,x_1\in D_1(f_p)$.
Say $p\forces \ff(x_0)=y_0$ and $\ff(x_1)=y_1$.
As $M_2\in\Khom$,
it suffices to find an approximation $h:[x_0,x_1]\ra [y_0,y_1]$
with $D_1(h)\contains D_1(f_q)$,
since then $k\circ h$ would be an approximation extending $g$ for
any order-preserving isomorphism $k:I_2\r (y_0,y_1)\ra I_2\r (g(x_0),g(x_1))$.

For $a\in (x_0,x_1)$, let $PV(a)$ denote the set of possible values of $\ff(a)$, 
i.e., the set of all $b\in (y_0,y_1)$ such that $r\forces \ff(a)=b$ for some
$r\ge p$.   If $a<a'$ then as $\forces \ff(a)<\ff(a')$, there always exist elements
$b\in PV(a)$
and $b'\in PV(a')$ such that $b<b'$.
By contrast, we say $PV(a)$ and $PV(a')$ {\it overlap\/} if there are
$c\in PV(a)$ and $c'\in PV(a')$ such that $c'\le c$.
Let the symmetric relation $R(a,a')$ hold if $PV(a)$ and $PV(a')$
overlap.
It is easy to verify that the set of elements $R$-related to $a$ is a convex subset
of $(x_0,x_1)$.
Let ${\sim}$ be the transitive closure of $R$.
For notation, let $[a]=\{a'\in (x_0,x_1):a\sim a'\}$.
Each $[a]$ is convex.
Similarly, for $b\in (y_0,y_1)$, let $PV(b)=\{a\in (x_0,x_1):f_p(a)=b\}$.
Define the relation $R$ on $(y_0,y_1)$ and $[b]$ analogously.
Note that if $b,c\in PV(a)$, then as $a\in PV(b)\cap PV(c)$, $[b]=[c]$.
As each of the equivalence classes are convex, it follows that 
if $R(a,a')$ holds and $b\in PV(a)$, $b'\in PV(a')$
then $[b]=[b']$.  Thus,  for all $a\in (x_0,x_1)$ and all
$b\in PV(a)$, $p\forces \ff:[a]\ra [b]$.
It is easy to see that if $a\in\dom(f_p)$ then $[a]=\{a\}$.
On the other hand,

\Claim{Subclaim}
If $a\not\in \dom(f_p)$ then there is a strictly increasing, cofinal and coinitial
sequence $\langle a_n:n\in\ZZ\rangle$ in $[a]$.

\bp
Suppose
$a\not\in\dom(f_p)$.  We first claim that there is an
$a'>a$, $a'\in [a]$.
To see this, pick  distinct elements $b_1,b_2\in PV(a)$
with $b_1<b_2$.  Pick $r\ge p$ with $r\forces \ff(a)=b_1$.
Pick $s\ge r$ with $b_2\in\dom(f_2^{-1})$ and let $a'=f_s^{-1}(b_2)$.
Then $a<a'$ and $b_2\in PV(a)\cap PV(a')$, so $a\sim a'$.
Similarly, there is $a'<a$ with $a'\in[a]$.
By symmetry, to complete the proof of the subclaim we need only show that
there is no strictly increasing
sequence $\langle a_\al:\al \in\o_1\rangle$ in $[a]$.
By way of contradiction, assume that such a sequence exists.
For each $\al\in\o_1$, let $A_\al=(x_0,a_\al)$ and let $B_\al=(y_0,b_\al)$.
As $\Q$ is countable, $PV(a')$ is countable for all $a'$, hence there is a club
$C\sbq \o_1$ such that, for all $\d\in C$, $a'\in A_\d$ implies 
$PV(a')\sbq B_\d$ and $b'\in B_\d$ implies $PV(b')\sbq A_\d$.
Thus, $p\forces \ff:A_\d\ra B_\d$ for $\d\in C$, 
contradicting the definition of $[a]$.
\endproof

Note that by symmetry, if $b\not \in \dom(f^{-1}_p)$
then there is a strictly increasing, cofinal and coinitial sequence of order type
$\ZZ$ in $[b]$.
We build our function $h:(x_0,x_1)\ra (y_0,y_1)$ as follows:
Let $h(a)=f_p(a)$ for each $a\in\dom(f_p)$.
For each non-trivial equivalence class $[a]$, let $b\in PV(a)$
and choose strictly increasing, cofinal and coinitial sequences
$\langle a_n:n\in\ZZ\rangle$ and $\langle b_n:n\in\ZZ\rangle$ in $[a]$ and $[b]$,
respectively.  Let $h(a_n)=b_n$ for each $n\in\ZZ$.
It is easy to verify that $h$ is an approximation on $[x_0,x_1]$.
\endproof

We close with the following example that shows that the assumption of
$D(T)$ countable in Theorem~\ref{big} cannot be weakened.

\begin{example} There is a  countable, superstable theory with a 
complete type of infinite multiplicity, yet non-isomorphism of
models of $T$ is preserved under all cardinal-preserving forcings.
\end{example}

Let $T$ be the theory of countably many binary splitting, cross-cutting
equivalence relations.  That is $L=\{E_n:n\in\o\}$ and the axioms of $T$
state that:
\begin{enumerate}
\item Each $E_n$ is an equivalence relation with two classes, each infinite and
\item For each $n\in\o$ and $w\sbq n$,
$\forall x\exists y (\bigwedge_{i\in w} E_i(x,y)\,\wedge\,\bigwedge_{i\in n\ss w}\neg
E_i(x,y)).$
\end{enumerate}

$T$ admits elimination of quantifiers, is superstable and the unique
1-type has infinite multiplicity.  However, for any model $M$ of $T$ and
any $a\in M$, every $p\in S_1(\{a\})$ is stationary.

Further, it is easy to verify that for all models $M,N$ of $T$ and all
$a\in M$, $b\in N$, there is an isomorphism $g:M\ra N$ with $g(a)=b$
if and only if for all 2-types $p(x,y)\in S_2(\emptyset)$,
$|p(M,a)|=|p(N,b)|$.

Now assume that $\forces_Q M\cong N$ for some cardinal-preserving forcing $Q$.
Then for some $q\in Q$, some $a\in M$ and some $b\in N$,
$$q\forces \text{``for all} p\in S_2(\emptyset),\  |p(M,a)|=|p(N,b)|.\hbox{''}$$

As $Q$ is cardinal preserving, this implies that $|p(M,a)|=|p(N,b)|$
for all $p\in S_2(\emptyset)$, so $M\cong N$.


\end{document}